\newcommand{\LTX}{\documentstyle{article} 
	   	  \newcommand{\Bbb}{\bf}
		  \newcommand{\myaddress}{paul@math.tau.ac.il} 
		 }

\newcommand{\AMSLTX}{\documentstyle[amstex]{article} 
		     \newcommand{\myaddress}{paul@math.tau.ac.il}  
		     }

\AMSLTX

\newcommand{\nline}{\ \\}

\newtheorem{thm}{Theorem}[section]

\newtheorem{lem}[thm]{Lemma}
\newtheorem{rem}[thm]{Remark}

\newtheorem{dfn}[thm]{Definition}
\newtheorem{cor}[thm]{Corollary}

\newcommand{\Qed}{\hfill \rule{.75em}{.75em}}

\newcommand{\cntrs}{\setcounter{thm}{0} \renewcommand{\thethm}
{\thesection.\Alph{thm}}}

\newcommand{\PoD}{Poincar\'{e} dual}
\newcommand{\Khlr}{K\"{a}hler}

\newcommand{\CPTU}{{{\Bbb C}P^2}}

\newcommand {\omstd}{\omega_{std}}
\newcommand {\Om}{\Omega}
\newcommand {\Ombar}{\overline{\Omega}}
\newcommand {\sig}{\sigma}
\newcommand {\Sig}{\Sigma}
\newcommand {\sigstd}{\sigma_{std}}
\newcommand {\lam}{\lambda}
\newcommand {\lambar}{\bar{\lambda}}
\newcommand {\Mbar}{\overline{M}}
\newcommand {\eps}{\epsilon}

\newcommand {\Th}{\Theta}
\newcommand {\vphi}{\varphi}

\newcommand {\cbar}{\bar{c}}
\newcommand {\Jbar}{\overline{J}}
\newcommand {\CC}{{\Bbb C}}
\newcommand {\NN}{{\Bbb N}}
\newcommand {\RR}{{\Bbb R}}
\newcommand {\QQ}{{\Bbb Q}}
\newcommand {\ZZ}{{\Bbb Z}}
\newcommand {\acs}{{almost complex structure}}

\newcommand {\NSWT}{{closed symplectic 4-manifold in the class ${\cal C}$}}
\newcommand {\bp}{b_2^+=1}
\newcommand {\calc}{{${\cal C}$}}

\newcommand {\MCDFROMLEMA}{lemma~2.2}

\title{Symplectic packing in dimension 4}
\author{Paul Biran\thanks{This is a part of my Ph.D. thesis, 
being carried out at Tel-Aviv University.}
\\ Tel-Aviv University \\ (\myaddress)}

\date{June 5, 1996}
      
\begin{document}

\maketitle

\begin{abstract}

We discuss closed symplectic 4-manifolds which admit full symplectic packings  
by $N$ equal balls for large $N$'s. We give a homological criterion for 
recognizing such manifolds. As a corollary we prove that ${\Bbb C}P^2$ can be 
fully packed by $N$ equal balls for every $N\geq 9$.

\end{abstract}

\section{Introduction} \cntrs

Let $(M^4,\Om)$ be a closed symplectic 4-manifold. We say that $(M,\Om)$ admits 
a symplectic packing by $N$ balls of radii $\lam_1,\ldots,\lam_N$ if there 
exists a symplectic embedding of the disjoint union 
$\coprod_{q=1}^N (B(\lam_q),\omstd)$ into $(M,\Om)$, where 
$(B(\lam_q),\omstd)$ denotes the standard closed 4-ball of radius $\lam_q$, 
endowed with the standard symplectic form of $\RR^4$, 
$\omstd=dx_1\wedge dy_1+dx_2\wedge dy_2$.

We say that $(M,\Om)$ admits a full packing by $N$ equal balls if the 
supremum of 
volumes which can be filled by symplectic embeddings of $N$ disjoint equal 
balls equals to the volume of $(M,\Om)$. Otherwise we say that there is a 
packing obstruction. Finally, we denote by $v_N(M,\Om)$ the ratio between the 
supremum of the fillable volume by packings with $N$ equal balls, and the 
volume of $(M,\Om)$.

Symplectic packings were studied for the first time by Gromov in~\cite{Gr}, 
and later by McDuff and Polterovich in~\cite{M-P}. 
McDuff and Polterovich discovered 
that for certain manifolds there are packing obstructions. Moreover, they
were able to compute $v_N$ of certain manifolds for many values of $N$. 
In particular they computed 
$v_N(\CPTU)$ for any $N$ which is square and for $1\leq N \leq 8$. It turned 
out that for every $N$ which is a square there exists full packing, while for 
every $1\leq N \leq 8$ which is not a square the are packing obstructions.
We refer the reader to~\cite{M-P} for more details about the symplectic 
packing problem.

In this paper we continue the above discussion, concentrating on manifolds 
which admit full packings by $N$ equal balls for large enough $N$'s. 
We give a homological condition for recognizing such manifolds, and a direct 
method for computing values of $N_0$, such that for all $N\geq N_0$ there 
exists full packings by $N$ equal balls. We then work out several examples 
including $\CPTU$, which turns out to admit full packings for every 
$N\geq 9$.  

The methods we use are based on the inflation procedure of 
Lalonde and McDuff and on Taubes theory of Gromov invariants via 
pseudo-holomorphic curves.

\section{Main results} \cntrs

Our main results are concerned with closed symplectic 4-manifolds of the 
following types:
\begin{itemize}
\item manifolds with $\bp$, $b_1=0$.
\item ruled manifolds and their blow-ups.
\end{itemize}

Manifolds of the above types belong to a wider class known as {\em manifolds 
which do not have SW-simple type}. We shall denote this class by \calc, 
remarking that our main results remain true for manifolds in this class. 

We refer the reader to~\cite{TAU-The_SW_and_Gromov,TAU-From-the-SW}, and 
to \cite{McD-Lectures,McD-From} for the precise definition of that class.
Finally, note that the class \calc\ is closed under the operation of 
blowing-up.

Given a symplectic manifold $(M,\Om)$, we shall denote by 
$c_1=c_1(TM,J)$ the first Chern class of the complex vector bundle $(TM,J)$, 
where $J$ is any \acs\ tamed by $\Om$.

\begin{dfn}
\label{dfn-set_D}

{\em 
Let $(M,\Om)$ be a closed symplectic 4-manifold. Consider the following set
$${\cal D}_{\Om} = \{B\in H_2(M;\ZZ) \mid \Om(B)>0,\, c_1(B)\geq 2,\, 
B\cdot B\geq 0\} .$$
Define 
$$d_{\Om}=\inf_{B\in {\cal D}_{\Om}} 
\frac{\Om(B)}{c_1(B)} \in [0,\infty] \,.$$
Here we use the convention that $\inf \emptyset = \infty$.
}

\end{dfn}

Before stating our main theorem we mention that given a symplectic manifold 
$(M,\Om)$, its volume is defined to be 
$Vol(M,\Om)=\int_{M} \frac{1}{2}\Om\wedge\Om$.

\begin{thm}
\label{thm-Main_1}

Let $(M,\Om)$ be a \NSWT. Suppose that $0<d_{\Om}\leq \infty$. Then 
$$ v_N(M,\Om) \geq \min \{1,\frac{Nd_{\Om}^2}{2Vol(M,\Om)}\}.$$
In particular, there exists an integer $N_{\Om}$ such that for every 
$N\geq N_{\Om}$, \  $(M,\Om)$ admits a full packing by $N$ equal balls. 
In fact, $N_{\Om}$ can be taken to be any integer which satisfies 
$$N_{\Om} \geq \frac{2Vol(M,\Om)}{d_{\Om}^2}.$$

\end{thm}

The proof is given in section~\ref{sect-Proof_main}, where a slightly 
sharper result is stated and proved.
The above theorem shows that it makes sense to define the following invariant.

\begin{dfn}
\label{dfn-d}

{\em 
Let $(M,\Om)$ be a closed symplectic 4-manifold. 
Define $P_{(M,\Om)}$, the {\em packing number} of $(M,\Om)$, as follows: 
$$P_{(M,\Om)} = 1 + \max \{N\in \NN \mid \mbox{there does not exist a full 
packing by $N$ equal balls}\}.$$ 
Here we use the convention that $\max \emptyset = 0$, while $\max$ of an 
unbounded set is $\infty$. When there is no risk of confusion we shall denote 
the packing number of $(M,\Om)$ by $P_{\Om}$.
}

\end{dfn}

As corollary to theorem~\ref{thm-Main_1} we prove:

\begin{cor}
\label{cor-Main_1}

1) $P_{(\CPTU,\sigstd)}=9$, where $\sigstd$ is the standard \Khlr\ form of 
$\CPTU$.\\
2) $2\frac{\beta}{\alpha}\leq 
P_{(S^2 \times S^2,\alpha \sig \oplus \beta \sig)} 
\leq 8 \frac{\beta}{\alpha}$, where $\sig$ is the standard symplectic form 
of $S^2$ and $0<\alpha\leq\beta$.\\
3) Let $R$ be an orientable surface of genus $g\geq 1$. 
Let $\Om=\beta \sig_R \oplus \alpha \sig_{S^2}$ be a symplectic form of 
$R\times S^2$, where $\sig_R$, $\sig_{S^2}$ are area forms of $R,S^2$ 
respectively such that $\int_R \sig_R = \int_{S^2} \sig_{S^2} = 1$.
Then $P_{\Om}=\lceil 2\frac{\beta}{\alpha}\rceil$. 
(Here, $\lceil x \rceil$ denotes the minimal integer which is grater or 
equal to $x$).

\end{cor}

The proof of this corollary is given in section~\ref{sect-Examples}.

\begin{rem} 
{\em
1) Notice that in part 3 of the above corollary, unlike in part 2, 
we do not assume that $\alpha\leq\beta$.\\
2) From parts 2 and 3 of the above corollary it follows that 
$P_{\Om}$ is not a deformation invariant, since $P_{\Om}\rightarrow\infty$ as 
$\alpha \rightarrow 0$ 
}
\end{rem}

\begin{thm}
\label{thm-Minimal_other}

Let $(M,\Om)$ be a closed symplectic minimal 4-manifold in the class \calc.
Suppose that $(M,\Om)$ is not rational or ruled, 
then ${\cal D}_{\Om}=\emptyset$. In particular $P_{(M,\Om)} = 1$.

\end{thm}

\begin{pf}
Denote by $K=-c_1(TM,J)$ the canonical class of $(TM,J)$, where $J$ is any 
\acs\ tamed by $\Om$. 
Since $(M,\Om)$ is minimal but not rational or ruled, then by a theorem 
due to Liu \cite{Liu-Some_new} we must have $K^2 \geq 0$ and 
$K\cdot[\Om] \geq 0$. Hence $K$ belongs to the closure of the 
positive light cone 
$$\overline{{\cal P^+}}=
\{a\in H^2(M;\RR)\mid a^2\geq 0,[\Om]\cdot a\geq 0\}.$$ 

Assume that ${\cal D}_{\Om}$ is not empty, and let $B\in {\cal D}_{\Om}$. 
Clearly ${\cal D}_{\Om} \subseteq \overline{{\cal P^+}}$. Since manifolds 
in the class \calc\ have $\bp$, it follows from the 
light cone lemma (see~\cite{McD-Lectures}) that $K \cdot B \geq 0$. But this 
cannot hold by the definition of the set ${\cal D}_{\Om}$. 
\ \Qed

\end{pf}

Examples of manifolds satisfying the conditions of the above theorem are:
Barlow surfaces, Dolgachev surfaces, hyper-elliptic surfaces and 
Enriques surface. 
As explained in~\cite{M-S-Survey} 
and~\cite{McD-From} the above manifolds belong to the class \calc. 
Since they are minimal and not rational or ruled they satisfy 
the conditions of theorem~\ref{thm-Minimal_other}.
\nline 

We conclude this section by mentioning that all symplectic packings of 
manifolds in the class \calc\ are unique in the sense 
that given $\lam_1,\ldots,\lam_N$, any two symplectic embeddings of the 
disjoint union of balls of radii $\lam_1,\ldots,\lam_N$ are symplectically 
isotopic. See~\cite{McD-Remarks} for the case of $\CPTU$ with 
$N\leq 2$, \cite{Bi-Connectedness} for the case of $\CPTU$ with $N\leq 6$,
and ~\cite{McD-From} for the general case.

\section{Inflation and Gromov invariants} \cntrs

In order to prove theorem~\ref{thm-Main_1} we need a few preliminary lemmas. 
For completeness we also state some relevant theorems from Taubes theory of 
Gromov invariants.

Let $(M,\Om)$ be a closed symplectic 4-manifold. 
An $\Om$-symplectic exceptional sphere is by definition a symplectically 
embedded sphere with self intersection number $-1$.
We shall denote by ${\cal E}$ the set 
of all 2-integral homology classes which can be represented by 
$\Om$-symplectic exceptional spheres.
In what follows we shall use the notation $PD$ for the {\PoD}ity. 

We first need the following theorem of McDuff which belongs to the framework 
of Taubes theory of Gromov invariants. For the theory of Gromov invariants 
see~\cite{TAU-The_SW_and_Gromov,TAU-From-the-SW},~\cite{McD-Lectures}.

\begin{thm}
\label{thm-Effective-class}

{\em (McDuff~\cite{McD-From} \MCDFROMLEMA)} Let $(M,\Om)$ 
be a \NSWT. Let $A\in H_2(M;\QQ)$ satisfy $\Om(A)>0$ 
and $A\cdot A > 0$. Then for sufficiently large $n$, the class 
$nA$ can be represented by a (possibly disconnected) 
$J$-holomorphic curve for generic $J$. 
If in addition, for every $E\in {\cal E}$ \  $A\cdot E \geq 0$ then there 
exists an embedded and connected pseudo-holomorphic curve representing 
the class $nA$.

\end{thm}

Combining this theorem with the inflation procedure of Lalonde and McDuff one 
obtains the following theorem: 

\begin{thm}
\label{thm-Strong_inflation}

{\em (McDuff~\cite{McD-From})} Let $(M,\Om)$ be a \NSWT. 
Let $A\in H_2(M;\QQ)$ satisfy 
$\Om(A)>0$ and $A\cdot A > 0$. Assume that for every 
$E\in {\cal E}$ \  $A\cdot E \geq 0$. Then there exists a closed 2-form 
$\rho$, representing $PD(A)$ and such that $\Om+y\rho$ is symplectic for all 
$y\geq 0$.

\end{thm}
For more details about this and about the inflation procedure 
see~\cite{McD-From} and~\cite{L-M-The_classification}.

\begin{rem}
{\em
An obvious conclusion of the above theorem is:\\
Under the assumptions of the above theorem, arbitrarily close to $PD(A)$ there 
exist cohomology classes which represent symplectic forms in the same 
deformation class as $\Om$.
}
\end{rem}

We shall also need the following well known lemma  
(see~\cite{McD-The_structure},~\cite{M-P}): 

\begin{lem}
\label{lem-Exc}

Let $(M,\Om)$ be a closed symplectic 4-manifold. Denote by ${\cal E}$ the set 
of all homology classes which can be represented by $\Om$-symplectic 
exceptional spheres. Then:\\
1) ${\cal E}$ depends only on the deformation class of $\Om$.\\
2) Let ${\cal J}(\Om)$ be the space of all $\Om$-tamed smooth {\acs}s of $M$. 
Then there exists a dense (actually, even residual) subset 
${\cal J}_{\cal E} \subseteq {\cal J}(\Om)$ such that for every 
$J\in {\cal J}_{\cal E}$ all classes in ${\cal E}$ admit $J$-holomorphic 
representatives which are connected, embedded and of genus 0.\\ 
3) If $E',E''$ are distinct classes in ${\cal E}$ then 
$E'\cdot E'' \geq 0$. 

\end{lem}

Essential to the symplectic packing problem is the symplectic blow-up 
operation (see~\cite{M-P},~\cite{M-S-Introduction}). We shall work in the 
following setting.  
Let $(M,J)$ be a 4-dimensional almost complex manifold with $J$ integrable 
near $x_1,\ldots,x_N \in M$. 
Let $(\Mbar,\Jbar)\stackrel{\Th}{\rightarrow} (M,J)$ 
be the complex blow-up of $(M,J)$ at $x_1,\ldots,x_N$. 
Denote by $\Sigma_q = \Th^{-1}(x_q)$ \ $q=1,\ldots,N$ the exceptional 
divisors, and by $E_q \in H_2(\Mbar;\ZZ)$ their homology classes. Finally, we 
set $e_q=PD(E_q)$.

Recall that a symplectic form taming an \acs\ $J$ is said to be $J$-standard 
near $x\in M$ if the pair $(\Om,J)$ is diffeomorphic to the standard pair 
$(\omstd,i)$ of $\RR^4$, near $x$.

The following lemma is an obvious generalization of 
proposition 2.1.C from \cite{M-P}.

\begin{lem}
\label{lem-Deformation}

Let $(M,\Om)$ be a closed symplectic 4-manifold. Let $J$ be an \acs\ tamed by 
$\Om$, which is integrable near $x_1,\ldots,x_N \in M$ and suppose that $\Om$ 
is $J$-standard near $x_1,\ldots x_N$. Let $\mu_1(0),\ldots,\mu_N(0)$ be 
positive numbers and let 
$$\vphi=\coprod_{q=1}^N \vphi_q : \coprod_{q=1}^N (B(\mu_q(0)),\omstd) 
\rightarrow (M,\Om)$$ be a symplectic embedding which is $(i,J)$-holomorphic.
Denote by $(\Mbar,\Ombar_0)$ the symplectic blow-up of $\Om$ with respect 
to $\vphi$. Suppose we have a symplectic deformation 
$\{\Ombar_t\}_{0\leq t\leq 1}$ starting with $\Ombar_0$, and lying in the 
cohomology class 
$$[\Ombar_t]=[\Th ^* \Om]-\pi \sum_{q=1}^N \mu_q(t)^2 e_q.$$ 
Then $(M,\Om)$ admits a symplectic packing by $N$ balls of radii 
$\mu_1(1),\ldots,\mu_N(1)$.

\end{lem}

The proof of this lemma goes along the same lines as those of proposition~2.1.C 
from~\cite{M-P}, only that here one has to adjust the symplectic forms 
$\Ombar_t$ on the exceptional divisors so that they become standard. 
This can be done using the symplectic neighborhood theorem.

\section{Proof of the main theorem} \cntrs
\label{sect-Proof_main}

We shall prove a slightly more general version of theorem~\ref{thm-Main_1}, 
namely:

\begin{thm}
\label{thm-Main_1a}

Let $(M,\Om)$ be a \NSWT. Suppose that $0<d_{\Om}\leq \infty$ and let  
$\lam_1,\ldots,\lam_N \, < \sqrt{d_{\Om}}$ be positive numbers which 
satisfy $$\sum_{q=1}^N \lam_q^4 < 2Vol(M,\Om).$$ Denote by 
$\Mbar \stackrel{\Th}{\rightarrow} M$ a complex blow-up of $M$ at 
$N$ distinct points. Then the following holds:\\
1) The cohomology class 
$$[\Th^*\Om]-\sum_{q=1}^N \lam_q^2 e_q$$ admits a symplectic 
representative.\\
2) The manifold $(M,\pi\Om)$ admits a 
symplectic packing by $N$ balls of radii $\lam_1,\ldots,\lam_N$.\\
3) In particular, if $$N\geq\frac{2Vol(M,\Om)}{d_{\Om}^2}$$ then there 
exists a full packing of $(M,\Om)$ by $N$ equal balls.

\end{thm}

\begin{pf} 
The idea of the proof goes along the following lines.
First, endow $\Mbar$ with some auxiliary symplectic form $\Ombar$, obtained 
from blowing-up $(M,\Om)$ with respect to some 
symplectic embedding of very small balls. Next, consider the homology class 
$A=PD([\Th^*\Om])-\sum_{q=1}^N \lam_q^2 E_q$. The idea is to show that for 
large enough $n$ the class $nA$ represents an embedded, reduced and connected 
pseudo-holomorphic curve in $\Mbar$. Once this is proved we can use the
inflation procedure to obtain a closed 2-form $\rho$ which lies in the 
cohomology class $[\Th^*\Om]-\sum_{q=1}^N \lam_q^2 e_q$ such that for every 
$y\geq 0$ the form $\Ombar+ y\rho$ is symplectic. Dividing by $y$ we see that 
for every $y>0$ the form $\frac{1}{y}\Ombar + \rho$ is symplectic. 
By taking $y$ to 
be very large we obtain a symplectic form lying in a cohomology class which is 
very close to the desired class: $[\Th^*\Om]-\sum_{q=1}^N \lam_q^2 e_q$.
It turns out that this approximation is enough for our purposes.
We now give the precise details of the proof.
\nline

Let $J$ be an \acs\ tamed by $\Om$ which is integrable near 
$x_1,\ldots,x_N \in M$. Let $(\Mbar,\Jbar) \stackrel{\Th}{\rightarrow} (M,J)$ 
be the complex blow-up of $M$ at $x_1,\ldots,x_N$. Denote by $\Sig_q$ the 
exceptional divisors, by $E_q$ their homology classes, and by $e_q$ the 
\PoD\ of $E_q$. Finally, denote by $c_1$ the first Chern class of $(TM,J)$ and 
by $\cbar_1$ the first Chern class of $(T\Mbar,\Jbar)$. Notice that 
$\cbar_1=c_1-\sum_{q=1}^N e_q$, under the natural decomposition 
$H^2(\Mbar;\ZZ)=H^2(M;\ZZ)\oplus \ZZ e_1 \oplus\cdots+\oplus \ZZ e_N$.

Without loss of generality we may assume that $\Om$ is $J$-standard near the 
points $x_1,\ldots,x_N$ since $\Om$ is isotopic to such a form 
(see~\cite{M-P} proposition~2.1.A).

The proof is divided into two steps. In the first we assume that the 
cohomology class of $\Om$ is rational. The second step is a reduction to the 
first one.

{\bf Step 1:} Assume that $[\Om]$ is a rational class.\\
Let $N\in \NN$, and let $\lam_q$ be as in the assumptions of the theorem.
Choose $\lambar_q \in \QQ$ such that $\lam_q<\lambar_q < \sqrt{d_{\Om}}$ and 
such that $\sum_{q=1}^N \lambar_q^4 < 2Vol(M,\Om)$. 
Set $$a = [\Th^*\Om]-\sum_{q=1}^N \lambar_q^2 e_q \in H^2(\Mbar;\QQ).$$
Denote by $A$ the \PoD\ of $a$. Clearly $A\cdot A > 0$.

Let $\Ombar_{\eps}$ be a symplectic form on $\Mbar$ obtained from symplectic 
blowing-up with respect to a holomorphic and symplectic embedding 
of $N$ balls of very small radii. Hence 
$$[\Ombar_{\eps}]=[\Th^*\Om]-\eps \sum_{q=1}^N e_q.$$

By taking $\eps$ be small enough we may assume that $\Ombar_{\eps}(A)>0$. 
Denote by ${\cal E}$ the set of homology classes which can be 
represented by $\Ombar_{\eps}$-symplectic exceptional spheres in $\Mbar$.

We claim that $A\cdot E>0$ for any $E\in {\cal E}$. 
For this purpose, write $$E=B-\sum_{q=1}^N m_q E_q \, , \, \mbox{ where } 
B\in H_2(M;\ZZ).$$

Notice first that $\Om(B)\geq 0$. Indeed, the form $\Ombar_{\eps}$ can be 
included in a smooth family of symplectic forms 
$\{\Ombar_t\}_{0< t \leq \eps}$ such that 
$$[\Ombar_t]=[\Th^*\Om]-t\sum_{q=1}^N e_q \,\,\,\,\, \mbox{for all}
\,\,\,0<t\leq \eps.$$ 
Since the set ${\cal E}$ depends only on the deformation class of 
$\Ombar_{\eps}$ we must have $\Ombar_t(E)>0$ for all $0<t\leq\eps$ and it 
easily follows that $\Om(B)\geq 0$.

Moreover, one can show that:\\
(i) If $\Om(B)=0$ then $B=0$ and $E=E_q$ for some $1\leq q \leq N$.\\
(ii) If $B\neq 0$ then $m_q\geq 0$ for all $1\leq q \leq N$.\\
Indeed, suppose that $\Om(B)=0$ but $B\neq 0$. Therefore $E\neq E_q$
for all $q$, so by lemma~\ref{lem-Exc} we have $E\cdot E_q \geq 0$,  
hence $m_q\geq 0$ for all $q$.
It is not possible that all the $m_q$'s are zero, since then $E=B$ and we get  
$\Ombar_{\eps}(E)=\Om(B)=0$. So there must be at least one positive $m_q$.
But then $\Ombar_{\eps}(E) = 0 - \eps \sum_{q=1}^N m_q < 0$, which is 
impossible. This proves that $B=0$, and it easily follows 
that $E=E_q$ for some $q$.
Part (ii) follows immediately from lemma~\ref{lem-Exc}.

Now we are ready to show that $A \cdot E > 0$.
If $\Om(B)=0$ then $E=E_q$ for some $q$, hence $A \cdot E = \lambar_q^2 > 0$. 
So assume $\Om(B)>0$, from which it follows that $m_q\geq 0$ for all $q$. 
If all the $m_q$ are zero then clearly 
$A\cdot E = \Om(B)=\Ombar_{\eps}(B) > 0$ 
and we are done. So assume that 
at least one of the $m_q$'s is positive, hence 
$$\sum_{q=1}^N m_q \geq 1 \, .$$

Since $\cbar_1(E)=1$ we must have $$c_1(B) = 1 + \sum_{q=1}^N m_q \geq 2.$$
Furthermore, since $E \cdot E = -1$ we must have 
$$B\cdot B - \sum_{q=1}^N m_q^2 = -1$$ hence $B \cdot B \geq 0$. 
This shows that $B \in {\cal D}_{\Om}$. 
Set $\lambar= \max \{\lambar_1, \ldots, \lambar_N\}$. We have 
$$A\cdot E = \Om(B)-\sum_{q=1}^N \lambar_q^2 m_q \geq
\Om(B)- \lambar^2 (c_1(B)-1) \geq 
\Om(B)-d_{\Om}c_1(B) + \lambar^2> 0.$$

Notice that the case 
${\cal D}_{\Om}=\emptyset$ (which corresponds to $d_{\Om}=\infty$), 
has already been treated since when ${\cal D}_{\Om}=\emptyset$, 
the proof of our claim ends at an earlier stage.

\nline

By theorem~\ref{thm-Strong_inflation} there exists a closed 2-form $\rho$,  
representing the class $a=PD(A)$ such that $\Ombar_{\eps} + y \rho$ is 
symplectic for all $y\geq 0$. As mentioned before the idea now is to consider 
the symplectic forms $\frac{1}{y}\Ombar_{\eps} + \rho$. 
If we take $y$ large enough we obtain in this way a symplectic form which 
lies in a cohomology class very close to the desired one. 

More precisely, consider the symplectic forms
\begin{equation}
\Ombar_s = \frac{1}{1+s-\eps}[ \Ombar_{\eps} + (s-\eps)\rho ] \; , \;
\;\;\;\;\;   \mbox{where \ } s \geq \eps .
\label{eq-Main_1a_1}
\end{equation} 

They lie in the cohomology class 
$$[\Ombar_s] = [\Th^*\Om] - \frac{1}{1+s-\eps}\sum_{q=1}^N 
(\eps+(s-\eps)\lambar_q^2) e_q .$$

Choose $s_0>\eps$ so large that 
$$\lam_q^2<\frac{\eps+(s_0-\eps)\lambar_q^2}{1+s_0-\eps}<\lambar_q^2 \;\;\;
\mbox{(we may assume $\eps<\lambar_q^2$ for all $q$)}.$$

Notice that equation~\ref{eq-Main_1a_1} provides a symplectic deformation 
$\{\Ombar_s\}_{\eps\leq s \leq s_0}$ starting at $\Ombar_{\eps}$ and ending 
with $\Ombar_{s_0}$, and this deformation satisfies the conditions of 
lemma~\ref{lem-Deformation}. 
It follows by this lemma that $(M,\pi\Om)$ admits a symplectic packing by 
$N$ balls of radii 
$\sqrt{ \frac{\eps+(s_0-\eps)\lambar_q^2}{1+s_0-\eps}} > \lam_q
\;\; (q=1,\ldots,N)$. 
In particular all the conclusions of the theorem hold for 
$\lam_1,\ldots,\lam_N$.

{\bf Step 2:} Consider the general case.\\
First notice that the proof of step 1 still holds if we assume 
that the cohomology class of $\Om$ is a real multiple of a rational class. 
We need now the following general and simple observation:\\
Let $a\in \RR^n$. Then arbitrarily close to $a$ there exist  
$a_0,\ldots,a_n \in \QQ^n\subset\RR^n$ and nonnegative real numbers 
$\alpha_0,\ldots,\alpha_n \geq 0$ such that $\sum_{i=0}^n \alpha_i = 1$ and 
$\sum_{i=0}^n \alpha_i a_i = a$.

Using this, we decompose 
$a=[\Th^*\Om]-\sum_{q=1}^N \lambar_q^2 e_q$ into a sum of the form 
$\sum_{i=0}^n\alpha_i a_i$ where $a_i \in H^2(\Mbar;\QQ)$, and $n=b_2(M)+N$. 
Since the classes $a_i$ can be 
taken to be arbitrarily close to $a$ we can proceed as we did in step 1 of the
proof. The only difference is that now we have to do $n+1$ inflations, thus 
obtaining $n+1$ closed 2-forms $\rho_0,\ldots,\rho_n$ representing the classes 
$a_0,\ldots,a_n$, and such that $\Ombar_{\eps}+\sum_{i=0}^N y_i\rho_i$ are 
symplectic for all $y_i \geq 0$. 
The proof now proceeds as in step 1 by taking $y_i=s \alpha_i$, 
with $s$ very large.

Note that inflation using more than one curve is slightly more complicated 
than just with one curve. We refer the reader to~\cite{McD-From} for more 
details on this.
  
Finally, the third statement of the theorem follows easily from the second. 

\ \Qed

\end{pf}

\section{Examples} \cntrs
\label{sect-Examples} 

In this section we work out some examples of computing the packing number.

\subsection*{Minimal rational and ruled manifolds} 

Consider first minimal rational manifolds, that is 
$\CPTU$ and $S^2\times S^2$. We start with $(\CPTU,\sigstd)$, where 
$\sigstd$ is the standard \Khlr\ form normalized such that 
$\int_{\CC P^2} \sigstd=1$.
Denote by $L$ the homology class of a projective line in $\CPTU$ and by $l$ 
its \PoD. We have $c_1=3l ,\: [\sigstd]=l$, hence $d_{\sigstd}=\frac{1}{3}$.
Using Theorem~\ref{thm-Main_1a} we see that given 
$\lam_1,\ldots,\lam_N < \frac{1}{\sqrt{3}}$ which satisfy 
$\sum_{q=1}^N\lam_q^4 < 1$, there exists a symplectic packing of 
$(\CPTU,\pi\sigstd)$ by $N$ balls of radii $\lam_1,\ldots,\lam_N$.

In particular, if $N\geq 9$, $\CPTU$ admits a full packing by $N$ equal balls.
Note however, that for 8 equal balls there is a packing obstruction 
(see~\cite{M-P}).
Finally, recall that by a theorem of Taubes any two cohomologous symplectic 
structures on $\CPTU$ are diffeomorphic. All the above prove: 

\begin{cor}
\label{cor-CPTU} 

For any symplectic form $\sig$ on $\CPTU$, $P_{\sig}=9$.

\end{cor}

Let us consider now $S^2\times S^2$. Again, it is enough to consider the 
standard split forms, since by theorems of McDuff and Li and Liu 
(see~\cite{L-M-The_classification}), any two 
cohomologous symplectic forms of $S^2 \times S^2$ are diffeomorphic.

Consider the symplectic form $\Om=\alpha \sig \oplus \beta \sig$, where 
$\int_{S^2}\sig = 1$ and $0<\alpha\leq\beta$. 
Using Gromov's non-squeezing theorem (see~\cite{Gr}), it is not hard 
to see that if $B(\lam)$ embeds symplectically into 
$(S^2\times S^2,\Om)$ then $\pi\lam^2 < \alpha$, hence we must have 
$$P_{\Om} \geq 2\frac{\beta}{\alpha}.$$
Using theorem~\ref{thm-Main_1a} we compute an upper bound for $P_{\Om}$ as 
follows:\\
Set $A_1=[S^2\times pt] ,\; A_2=[pt\times S^2]$, and $a_i=PD(A_i) \;i=1,2$. 
We have $c_1=2(a_1+a_2)$ and $Vol(S^2\times S^2, \Om)=\alpha\beta$.
Let $B\in {\cal D}_{\Om}$, say $B=n_1A_1+n_2A_2$. It is easy to see that 
$n_1,n_2$ are non-negative. 
Hence $$\frac{\Om(B)}{c_1(B)}=\frac{\alpha n_1+\beta n_2}{2n_1+2n_2}\geq
\frac{\alpha}{2}>0.$$ Therefore $d_{\Om}\geq \frac{\alpha}{2}$, hence 
from theorem~\ref{thm-Main_1a} we obtain:

\begin{cor}
\label{cor-S2xS2}

Let $0<\alpha\leq\beta$ and let $\sig$ be any area form of $S^2$. 
Then $2\frac{\beta}{\alpha} \leq P_{\Om} \leq 8\frac{\beta}{\alpha}$.

\end{cor}

Note that in some cases it is possible to give sharper bounds for $P_{\Om}$, 
and even to compute its precise value.
It is well known that there exists a diffeomorphism between 
the blow-up of $S^2\times S^2$ at $N$ points and the blow-up of $\CPTU$ at 
$N+1$ points. Furthermore it is not hard to compute that this 
diffeomorphism can be chosen to induce the following correspondence:\\ 
packing $(S^2 \times S^2, \pi\Om)$ by $N$ balls of 
radii $\lam_1,\ldots,\lam_N$ correspond to packing 
$(\CPTU,\pi(\alpha+\beta-\lam_1^2)\sigstd)$ by $N+1$ balls of radii 
$\sqrt{\alpha-\lam_1^2}, \sqrt{\beta-\lam_1^2}, \lam_2,\ldots,\lam_N$.

Restricting to $\alpha=\beta$, a straightforward computation shows that there 
is packing obstruction for 7 balls (one uses the above correspondence and the 
packing inequalities from~\cite{M-P}). Hence, if $\alpha=\beta$ 
then $P_{\Om}=8$.\\

Let us consider now irrational ruled surfaces.

\begin{cor}
\label{cor-RxS2}

Let $R$ be an orientable surface of genus $g\geq 1$. 
Let $\sig_R$, $\sig_{S^2}$ be area forms on $R$,$S^2$ respectively, such that 
$\int_R \sig_R = \int_{S^2} \sig_{S^2} = 1$. Let $\alpha,\beta$ be positive 
numbers, then 
$$P_{(R\times S^2,\beta \sig_R \oplus \alpha \sig_{S^2})}=
\lceil 2\frac{\beta}{\alpha} \rceil$$

\end{cor}

\begin{pf}
Set $\Om=\beta \sig_R \oplus \alpha \sig_{S^2}$. 
Let $J_s$ be a split complex structure on $R\times S^2$. 
Denote by $((R\times S^2)_N, \Jbar_s)$ the complex blow-up of $R\times S^2$ 
at $N$ generic points, and by $\Ombar$ a blow-up of $\Om$. 
We claim that the set ${\cal E}$ of $\Ombar$-symplectic exceptional spheres 
is: 
$${\cal E}=\{E_1,\ldots,E_N,S-E_1,\ldots,S-E_N\},$$ where $S=[pt\times S^2]$.

Indeed, let $C$ be an $\Ombar$-symplectic exceptional sphere in the class $E$.
Choose a generic \acs\ $\Jbar$ tamed by $\Ombar$ for which $C$ is 
$\Jbar$-holomorphic. 

Consider a generic path $\{\Jbar_t\}_{0\leq t \leq 1}$ of $\Ombar$ tamed 
{\acs}s with $\Jbar_0=\Jbar$ and $\Jbar_1=\Jbar_s$.
Since $\Jbar$ can be chosen to be arbitrarily close to $\Jbar_s$,  
we may assume that for all $0\leq t<1$ there exist $\Jbar_t$-holomorphic 
$E$-spheres. Using Gromov's compactness theorem we obtain a
(possibly cusp) $\Jbar_s$-holomorphic $E$ curve, 
$\widetilde{C}=C_1\cup \ldots \cup C_n$, with $genus(C_i)=0$.

Denote by $pr_{_R}:R\times S^2 \rightarrow R$ the projection and by 
$\widetilde{pr_{_R}}$ its lifting to $(R\times S^2)_N$. 
Clearly $\widetilde{pr_{_R}}$ is $\Jbar_s$-holomorphic.
Since $genus(R)\geq 1$ it follows that for every $j$, \ 
$\widetilde{pr_{_R}}(C_j)$ 
must be a point, say $p_j$. As $\widetilde{C}$ is connected we see that 
all the $p_j$'s are equal. Thus $\widetilde{pr_{_R}}(\widetilde{C})$ 
is a point, and it immediately follows that 
$E\in\{E_1,\ldots,E_N,S-E_1,\ldots,S-E_N\}$. Conversely, it is obvious that 
$E_q$ and $S-E_q$ are $\Ombar$-symplectic exceptional classes.\\

Since the exceptional spheres provide all the packing obstructions, 
it follows that $(R\times S^2,\pi\Om)$, admits a symplectic packing by $N$ 
balls of radii $\lam_1,\ldots,\lam_N$ iff $\lam_q^2 < \alpha$ and 
$\sum_{q=1}^N \lam_q^4 < 2\alpha\beta$. Restricting to equal balls the result 
easily follows.

\ \Qed

\end{pf}

\nline
\noindent {\em Acknowledgments.} 
I wish to thank Prof. Dusa McDuff for giving me the 
opportunity to read early versions of her 
papers~\cite{McD-From,McD-Lectures} which initiated the work on this paper,
for helping me understand the relevant parts of Taubes theory,
and for very useful remarks. 
I am also grateful to my teacher Prof. Leonid Polterovich for his 
continual encouragement and for many hours of extremely useful conversations.
Finally, I would like to thank Roman Bezrukavnikov for simplifying the 
argument of the proof of corollary~\ref{cor-RxS2}.

\end{document}